\documentclass[11pt,a4paper]{article}
\usepackage{makeidx}
\usepackage{graphicx}
\usepackage{latexsym}
\usepackage{amsmath}
 \usepackage{amsfonts}
\usepackage{verbatim}
\usepackage{nameref}
\usepackage{hyperref}
\usepackage{sectsty}
\usepackage{enumerate}
\usepackage{subfigure}
\usepackage{epsf}
\usepackage{amssymb}

\usepackage{amsmath}
 \usepackage[dvips]{epsfig}
 \usepackage{url}
\usepackage{slashbox}

 \usepackage{float,caption}
\newcommand\blam{\mbox{\boldmath$\lambda$}}

\newcommand{\real}{{\bf R}}

\newcommand{\half}{\frac{1}{2}}

\newcommand{\bsig}{\mbox{\boldmath$\sigma$}}

\newcommand{\barsig}{\bar{\sigma}}

\newcommand{\bM}{{\bf M}}
\newcommand{\bK}{{\bf K}}
\newcommand{\eb}{\begin{equation}}
\newcommand{\ee}{\end{equation}}

\newcommand{\calP}{{\cal{P}}}

\newcommand{\bG}{{\bf G}}
\newcommand{\bH}{{\bf H}}

\newcommand{\bff}{{\bf f}}

\newcommand{\bw}{{\bf w}}

\newcommand{\calS}{{\cal S}}

\newcommand{\calU}{{\cal U}}

\newcommand{\calR}{{\cal R}}

\newcommand{\barw}{{\bar{w}}}

\newcommand{\alp}{{\alpha}}

\newcommand{ \eps}{{\epsilon}}

\newcommand{ \sig}{{\sigma}}
\newcommand{ \Lam}{{\Lambda}}

\newcommand{ \lam}{{\lambda}}

\newtheorem{theorem}{Theorem}

\renewcommand\real{ \mathbb{R}}
\renewcommand\calR{\real}
\newcommand\barbw{\bar{\bw}}
\newcommand\barbsig{\bar{\bsig}}
\begin{document}
\vspace{1cm}
\begin{center}
{\Large {\bf On SDP Method for Solving Canonical Dual  Problem in Post Buckling  
of Large Deformed  Elastic Beam  \vspace{0.5cm}\\}}
\vspace{0.5cm}

{\bf Elaf Jaafar Ali$^a$} \&
{\bf David  Yang Gao$^{b}$}\\

\vspace{0.5cm}

{\em  $^{a,b}$Faculty of Science and Technology,\\ Federation University Australia, Mt Helen, Victoria 3353, Australia}\\
{\em  $^{a}$University of Basrah, College of Science, Basra, Iraq}\\
{\em  $^a$elafali@students.federation.edu.au}  \;\;\&\;\;
{\em $^{b}$d.gao@federation.edu.au}\\ \vspace{.5cm}
\end{center}

\begin{abstract}   
   This paper presents a new methodology and algorithm for solving post buckling problems of a large deformed elastic  beam.  The   total potential energy of  this beam  is a nonconvex functional, which can be used  to model  both  pre- and post-buckling phenomena.   By using a canonical dual finite element method,  a new primal-dual semi-definite programming  (PD-SDP) algorithm is presented, which can be used to obtain all possible post-buckled solutions. Applications are illustrated by several numerical examples  with different boundary conditions. We find that the global minimum solution of the nonconvex potential leads to a stable configuration of the buckled beam, the local maximum solution leads to the unbuckled state, and both of these two solutions are numerically stable. However, the  local minimum solution leads to  an unstable buckled state, which  is  very sensitive to axial compressive forces,  thickness of beam,  numerical precision,   and the size of  finite elements. The method and algorithm proposed in this paper can be used for solving general nonconvex variational problems in engineering and sciences.      
\end{abstract}


\section{Introduction} \label{sec:Introduction}
It is known that the total potential energy for the post-buckling of large deformed structures must  be nonconvex to allow multiple
local minimum solutions for all possible  buckled status \cite{gao-mrc96}.
 However, nonconvex   variational problems have always presented  serious challenges  not only
  in computational mechanics, but also in  mathematical analysis   and computer science \cite{gao-amma16}.
  Traditional finite element methods for solving any nonconvex variational problem  usually end up with
   a nonconvex minimization problem in $\mathbb{R}^n$.
  Due to the lack of global optimality criteria,
   popular nonlinear programming methods developed from  convex optimization  can't be used  to find  global optimal solutions.
   It was discovered in \cite{gao-ogden-qjmam} that for certain  external loads,
    both global and local minimum  solutions to large deformed mechanics problems are usually nonsmooth and can't be captured by
  any Newton-type methods.
  Therefore, most nonconvex  optimization problems are considered as {\em NP-hard} (Non-deterministic Polynomial-time hard) in computer science.
Unfortunately, these well-known difficulties are not fully recognized  in computational mechanics due to the significant  gap between
  engineering mechanics and global optimization. Indeed,
  engineers and scientists are mistakenly attempting to use traditional finite element methods and commercial softwares
 for solving nonconvex mechanics problems.

Canonical duality theory is  a newly developed and potentially powerful methodology which can be used  not only for modeling complex systems within a unified framework, but also for solving a large
class of challenging  problems in nonconvex, nonsmooth, and discrete  systems \cite{gao-amma16}.  This theory  comprises mainly three parts:
1)  a canonical dual transformation,  which can be used to formulate perfect dual problem without a duality gap;
2)  a complementary-dual variational principle, which presents a unified analytic solution form for general problems in continuous and discrete systems;
 3) a triality theory, which can be used to identify   both global and local extrema and to develop effective algorithms for solving nonconvex optimization problems.

  The   canonical duality theory was developed from   Gao and Strang's original work on
 nonconvex/nonsmooth variational/boundary value problems in finite deformation systems \cite{gao-strang}.
 In order to recover the  complementary energy principle  in  nonconvex analysis, they
discovered a so-called {\em complementary gap function}, which leads to a complementary-dual
variational principle in finite deformation mechanics. They proved that the positivity of  this gap function
 provides a global optimality condition for nonconvex variational problem.
  It was realized  seven years later that
the negativity of this  gap function can be used to identify the biggest local minimal and local maximal solutions.
 Therefore, a triality theory was first proposed in post-buckling problems of a large deformation beam model \cite{gao-amr97}, and a
pure complementary energy principle was  obtained in 1999  \cite{gao-mrc99}.
This principle solved an open problem in nonlinear elasticity  \cite{li-gupta}, which can be
used for obtaining analytical solutions to general  large deformation problems \cite{gao-mecc,gao-anti,gao-amma16}.
Based on the canonical duality theory and the mixed finite element method, a canonical dual finite element method has been developed
\cite{gao-em96} with the  successful  application  for solving  nonconvex mechanics problems  in phase transitions of solids
\cite{gao-yu}.
 It was discovered recently   \cite{cai-gao-qin,santos-gao} that
the nonconvex variational problem of a   post-buckled nonlinear Gao beam
can have at most three smooth solutions: a global minimizer representing  a stable buckled state,
  a local maximizer for an unbuckled state, and a local minimizer for an unstable buckled state.
  Both global minimum and local maximum  solutions can be obtained easily by using the canonical dual finite element method.
  However,  the local minimum solution is very sensitive and  difficult to obtain by standard convex minimization algorithms.

  The main goal of the present paper is to develop a new canonical primal-dual algorithm for solving the post-buckling problem with  special attention to the local unstable buckled configuration of a large deformed beam.
 The generalized total complementary energy associated with this model is a nonconvex functional and is reformulated as a global optimization problem to study 
 the post-buckling responses of the beams.
  Based on the canonical duality theory and the associated triality theorem, a new primal-dual semi-definite program (PD-SDP) algorithm is proposed for solving this challenging problem  to obtain all possible
   solutions.
   Applications are illustrated by different boundary value problems. An important mistake in \cite{cai-gao-qin} on the local  minimum solution is found.

\section{Nonconvex  problem and canonical duality theory}
\label{sec:Gao Beam model}
\begin{figure}[htbp]
 \centering
\scalebox{0.14}
{\includegraphics{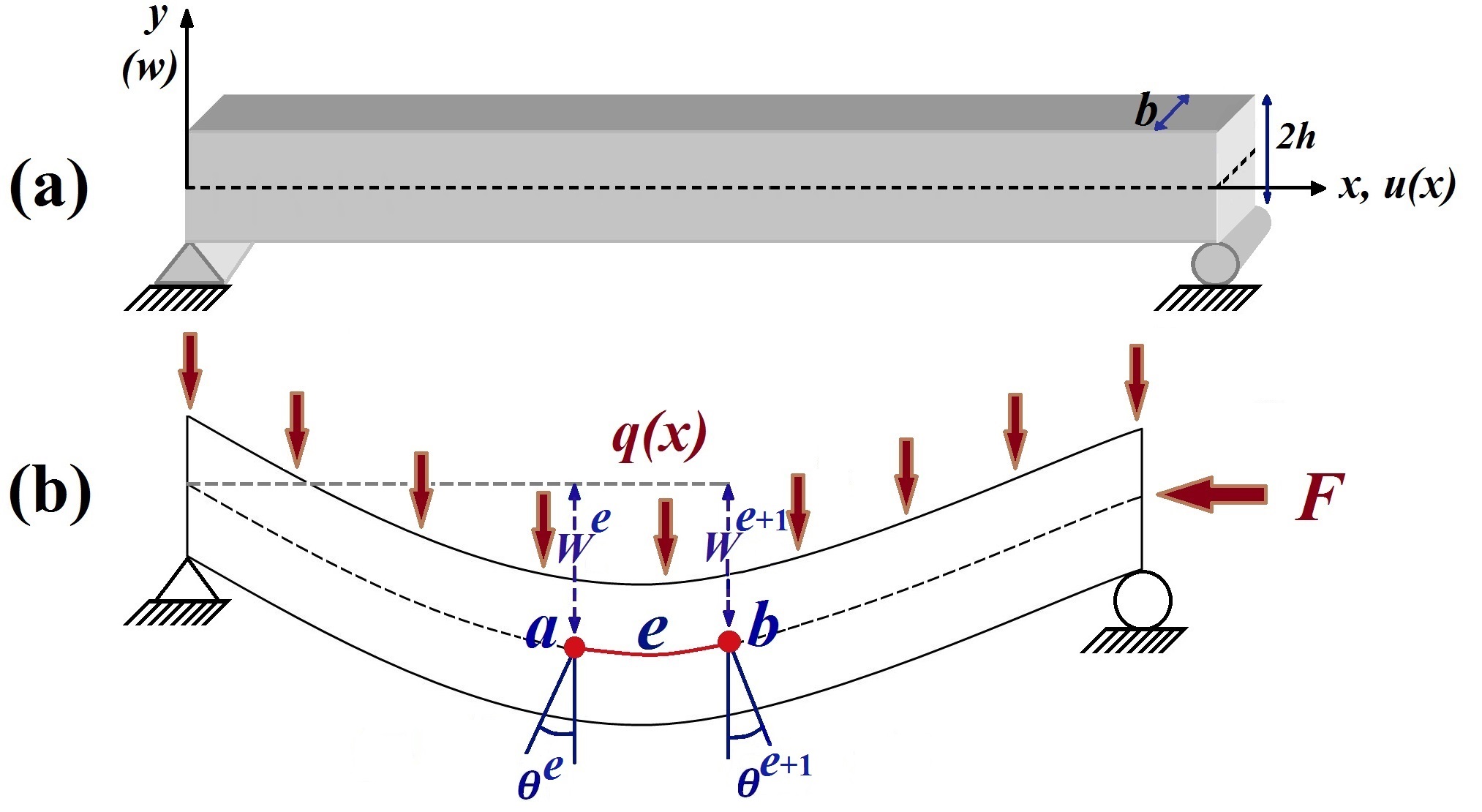}}
\caption{{\em Simply supported beam model - pre and post buckling analysis}}
\label{SSBeam}
\end{figure}

Let us consider an elastic  beam subjected to a vertical distributed lateral load $q(x)$ and compressive external axial force $F$ at the right end as shown in Figure  \ref{SSBeam}.
 It was discovered by Gao  in 1996 \cite{gao-mrc96} that the well-known von Karman nonlinear plate model in one-dimension
  is equivalent to a linear differential equation,  therefore, it can't be used to study post-buckling phenomena.
  The main reason for this ``paradox'' is due to the fact that the stress in the lateral direction of a large deformed plate was  ignored by von Karman.
  Therefore,  the von Karman equation works only for thin-plates and can't be used as a beam model.
  For a relatively thick beam  such that $h/L \sim w(x) \in O(1)$,
 the deformation in the lateral direction can't be ignored. Based on the finite deformation theory for Hooke's  material and the
 Euler-Bernoulli hypothesis (i.e.  straight lines normal to the mid-surface remain straight and normal to the mid-surface  after deformation),
a nonlinear  beam model was proposed by Gao   \cite{gao-mrc96}:
\eb\label{eq:one}
 E I w_{,xxxx}- \alpha E w^2 _{,x} w_{,xx}+E \lambda w_{,xx} -f(x)=0,  \\\   \forall x\in [0  ,    L]
\ee
where $E$ is the elastic modulus of material, $I = 2h^{3}/3$ is the second moment of area of the beam's cross-section,
 $w$ is the transverse displacement field of the beam,
$ \alpha= 3h(1 - \mu^{2}) > 0$  with $\mu$ as  the Poisson's ratio, 
$ \lambda = (1 + \mu)(1 - \mu^{2}) F/E > 0$ is
an  integral constant,
 $ f (x) = (1 - \mu^{2})   q(x) $ depends mainly on the distributed lateral load $q(x)$,
  $2h$ and $L$ represent  the height and length of the beam, respectively.
 The axial displacement $u(x)$   is governed  by the following differential equation \cite {gao-mrc96}:
\eb \label{eq-u}
u_{,x}=-\frac{1}{2}(1+\mu) w^{2}_{,x}-\frac{\lambda}{2h(1+\mu)}\; ,
\ee
which shows that if $u(x)  \sim w_{,x}(x) \in O(\eps), $ then $\; u_{,x}(x)  \sim w_{,xx}(x) \in O(\eps^2)$.

The total potential energy   of this beam model  is   $\Pi(w): \calU_{a}\rightarrow \calR $ defined by
\eb
 \Pi(w)=\int^{L}_{0}\bigg( \frac{1}{2} EIw^{2}_{,xx}+\frac{1}{12}E \alpha  w^4 _{,x} - \frac{1}{2}E \lambda w^{2}_{,x} -f(x) \; w\bigg) dx ,
\ee
where $\calU_{a}$ is the kinematically admissible space, in which  certain necessary boundary
conditions  are given.
Thus, for the given  external load $f(x)$ and end load $\lambda$, the primal variational problem  is  to find $\bar{w}\in \calU_{a}$ such that
\eb
   (\calP): \;\;\;\;\Pi(\bar{w})= \inf{\{ \Pi(w)|w\in \calU_{a}\}}.
\ee
It is easy to prove that the stationary condition
$ \delta \Pi( {w})=  0 $ leads to the governing equation (\ref{eq:one}).

If the nonlinear term in (\ref{eq:one}) is ignored and $f = 0$, then this nonlinear Gao beam is degeneralized to the well-known Euler-Bernoulli beam equation\footnote{Strictly speaking, instead of $\lambda$, the axial load in the Euler-Bernoulli beam should be  
$F = \lambda E/[(1+\mu)(1-\mu^2)]$.}: 
\eb\label{eq:Fcr2}
   EIw_{,xxxx}+ \lambda  E  w_{,xx}= 0 .
\ee
 It is known that this linear  beam will be buckled if  the axial load $\lambda$ reaches the Euler  buckling load $\lambda_{cr}$ defined by 
\eb\label{eq:Fcr}
   \lambda_{cr}= \inf_{w\in \calU_{a}}\frac{\int^{L}_{0} EIw^{2}_{,xx}dx}{\int^{L}_{0} E w^{2}_{,x}dx}. 
\ee

Clearly,     in the pre-buckling state, i.e. before the axial load $\lambda$ reaches  the Euler buckling load $\lambda_{cr}$,
we have
\eb
 \Pi_{EB}(w) = \int^{L}_{0} EIw^{2}_{,xx}dx   -  \lambda \int^{L}_{0} E w^{2}_{,x}dx > 0  \;\; \forall w \in \calU_a, \;\; \lambda < \lambda_{cr} .
\ee
 In this case,   $\Pi_{EB}(w)$ and 
 $\Pi(w)$ are  strictly  are strictly convex  on $\calU_{a}$, therefore, both the Euler-Bernoulli beam   (\ref{eq:Fcr2}) 
  and the nonlinear  Gao beam  (\ref {eq:one})  can have only one solution (see Lemma 2.1.
and Theorem 2.1 in \cite{m-n-17}).  
 
Dually, in the post-buckling state, i.e.  $ \lambda > \lambda_{cr}$,   the total potential energy for the Euler-Bernoulli beam is strictly concave and 
 \eb
\inf \left\{ \Pi_{EB}(w)  | \;\; w \in \calU_a, \;\; \lambda >  \lambda_{cr}   \right\} = - \infty,
\ee
which   means that the Euler-Bernoulli beam  is crushed. This shows that the  Euler-Bernoulli beam can't be used for studying post-buckling problems.
However,   for the nonlinear Gao beam, it was proved recently by Machalov\'{a}  and Netuka (see Remark 2.2,  \cite{m-n-17}) 
 that there exists a constant $\lambda_{cr}^G \ge \lambda_{cr}$ such that 
 the total potential energy $\Pi(w)$  is a  nonconvex (double-well) functional if $ \lambda > \lambda_{cr}^G$,
which allows at most  three critical points, i.e. the strong solutions to the
nonlinear equation (\ref{eq:one}) at  each material point $x\in [0 \   L]$:
 two minimizers corresponding to the two possible  buckled states,
  one local maximizer corresponding to the possible unbuckled state \cite{Gao-2008}.  
    Clearly, these  solutions  are  sensitive to both the axial load $\lambda$ and the distributed lateral force field $f(x)$.
By  equation (\ref{eq-u}) we know that the axial deformation could be relatively large,
while the nonconvexity of the total potential shows that  this nonlinear beam model can be used for studying both pre and post-buckling  problems \cite{cai-gao-qin,santos-gao}.
Recently, the Gao beam model has been generalized  for many real-world applications
  in engineering and sciences \cite{ahn-etal,and-etal,b-d-s17,k-l-s,leve,m-n,m-n-152,m-n-171}. 
 
 Although the nonlinear Gao beam can be used for modeling natural phenomena,  the  nonconvexity of this beam model leads to some fundamental challenges in mathematics and computational science. 
Generally speaking,
traditional numerical methods and  nonlinear optimization techniques can be used only for solving convex  minimization problems.
Due to the lack of a global optimality criterion to identify a global minimizer at each iteration,
  most  nonconvex optimization problems can't be solved deterministically, therefore,
  they are considered to be  NP-hard in global optimization and computer
  science \cite {Gao-Sherali}.

It was shown in  \cite {gao-book00} that by introducing a canonical strain measure $\epsilon = \Lambda(w) =  \half w^2_{,x}$
and a convex canonical function $V(\epsilon) = \frac{1}{3} E\alpha \epsilon^2 - E\lambda \epsilon$,
the nonconvex (double-well) potential $W(w_{,x}) = \frac{1}{12} E \alp w^4_{,x} -  \half E \lambda w^2_{,x}$  in $\Pi$
can be written in the
canonical form $W( w_{,x}) = V(\Lambda(w ))$. Thus, the canonical
 dual stress can be uniquely defined by
\eb\label{eq:sigma}
   \sigma= \partial V(\epsilon) = \frac{2 E\alpha}{3} \epsilon  - E\lambda .
\ee
By the Legendre transformation, we  have the canonical complementary energy
 \[
 V^*(\sig) = \epsilon \sigma - V(\epsilon) =    \frac{3}{4E\alpha}(\sigma+ E\lambda)^2 .
 \]
 Thus, replacing $W(w_{,x})$ with $V(\Lambda(w)) = \Lambda(w) \sig - V^*(\sig)$,
  the Gao-Strang total complementary energy $\Xi:\calU_{a}\times  \calS_{a}\rightarrow\calR$ \cite{gao-strang} in nonlinear elasticity can be defined as
\begin{eqnarray}
\Xi (w,\sigma) &=& \int_{0}^{L}\bigg(\frac{1}{2} EI w_{,xx}^{2} +\frac{1}{2} \sigma w^{2} _{,x} - \frac{3}{4E\alpha}(\sigma+ E\lambda)^2 -f(x) w\bigg) dx \nonumber \\
&= &  G(w,\sig) - \int_0^L [ V^*(\sig)- f(x) w] dx,\label{eq:XiXi}
\end{eqnarray}
where $\calS_{a} = \{ \sig \in C[0,L] | \;\; \sig(x) \ge - \lambda E \;\; \forall x \in [0,L]\}$
  and
\eb
G(w,\sigma)=\int_{0}^{L}\bigg(\frac{1}{2} EI w_{,xx}^{2} +\frac{1}{2} \sigma w^{2} _{,x}\bigg) dx 
\ee
is the  generalized Gao-Strang complementary gap function  \cite{gao-strang}.
 \begin{theorem}[Complementary-duality Principle]
  For any given external load $f(x)$ and end load $\lambda$, the pair  $(\barw, \barsig)$ is a critical point of $\Xi(w, \sig)$
 if and only if $\barw$ is a critical point of $\Pi (w)$ and $\Pi(\barw) = \Xi(\barw, \barsig)$.\vspace{-.3cm}
 \end{theorem}
 {\bf Proof}. The criticality condition $\delta \Xi(\barw, \barsig) = 0$ leads to the following canonical equations:
 \begin{eqnarray}
& &  E I \barw_{,xxxx} - \barsig \barw_{,xx} = f(x) , \;\; \label{eq-caneqns} \\
& & \half \barw_{,x}^2 = \frac{3}{2 E\alp} (\barsig + E \lam),\label{eq-cancons}
 \end{eqnarray}
 which are equivalent to  equation (\ref{eq:one}). The equality $\Pi(\barw) = \Xi(\barw, \barsig)$
 follows directly from the Fenchel-Young equality $V(\Lam(\barw) ) + V^*(\barsig) = \Lam(\barw) \barsig$
 due to the convexity of the canonical function $V(\epsilon)$.  \hfill $\Box$

 \begin{theorem}[Triality Theory] Let  $(\barw, \barsig)$ be  a critical point of $\Xi(w, \sig)$.

 If $G(\barw, \barsig) \ge 0  $, then $\barw$ is a global minimizer of $\Pi (w)$ on $\calU_a$ and
 \eb \label{eq-minmax}
 \Pi(\barw) = \min_{w \in \calU_a} \Pi(w) = \min_{w \in \calU_a } \max_{\sig \in \calS_a} \Xi(w, \sig). \vspace{-.3cm}
 \ee

  If $G(\barw, \barsig) <  0   $, then on a neighborhood $\calU_o \times \calS_o$ of $(\barw, \barsig)$,
   we have either
 \eb \label{eq-minx}
 \Pi(\barw) = \min_{w \in \calU_o} \Pi(w) = \min_{w \in \calU_o} \max_{\sig \in \calS_o} \Xi(w, \sig) = \Xi(\barw, \barsig)  ,  \vspace{-.3cm}  \ee
 or
 \eb \label{eq-maxx}
 \Pi(\barw) = \max_{w \in \calU_o} \Pi(w) = \max_{w \in \calU_o} \max_{\sig \in \calS_o} \Xi(w, \sig) = \Xi(\barw, \barsig) \vspace{-.3cm} .
 \ee
\end{theorem}
{\bf Proof}. For the positive gap function, $\Xi(w, \sig)$ is a saddle functional  and the total potential $\Pi(w)$ is convex on $\calU_a$ \cite{gao-strang}. In this case,  statement (\ref{eq-minmax}) follows directly from
Gao and Strang's  theory  for  general large deformation problems \cite{gao-strang}.
While for the negative gap function,  $\Xi(w, \sig)$ is a bi-concave functional. In this case, the total potential $\Pi(w)$ is nonconvex on $\calU_a$, which could have both local minimum  and local maximum solutions.
Due to the fact that  $  \max_{\sig \in \calS_o} \Xi(w,\sig) =   \max_{w \in \calU_o}\Pi(w)$,
the statements (\ref{eq-minx}) and  (\ref{eq-maxx}) can be proved easily by the general triality theory \cite{gao-book00}. \hfill $\Box$

The  {\em triality theory} was first discovered  in the post-buckling analysis of the large deformed beam theory \cite {gao-mrc99}.
 Generalization  to nonconvex/discrete optimization problems  was  given   in 2000\cite{gao-cace09}.
Detailed information relating to this theory as well as
 its extensive applications in global optimization as well as nonconvex mechanics
 can be found in  the monograph \cite {gao-book00} and recent review articles \cite{gao-cace09,gao-amma16,gao-bridge}.

\section{Mixed finite element method}
By using  the finite element method, the domain of the  beam is discretized into $m$ elements $ [0, L] = \bigcup_{e=1}^m \Omega^e$.
  In each element  $\Omega^e = [x^e,x^{e+1}]$, the deflection, rotating angular and dual stress for  node $x^e$ are marked as   $w^e$, $\theta^e$
and $\sigma^e$, respectively,  and similar for  node $x^{e+1}$.
Then, we have the nodal displacement vector $w_e$ of the e-th element 
\eb
w^{T}_{e}= [w^e \;\;  \theta^e  \;\;  w^{e+1} \; \;  \theta^{e+1}],
\ee 
 and the nodal dual stress element $\sigma_e$
\eb
\sigma^{T}_{e}= [\sigma^e \;\; \sigma^{e+1}].
\ee
 In each element, we use mixed finite element interpolations for both
 $w(x)$ and $\sig(x)$, i.e.
\[
w^{h}_{e}(x)= N^{T}_{w}(x)  w_{e} \;\;\;, \;\;\;\;\;   \sigma^{h}_{e}(x)=N^{T}_{\sigma}(x)  \textbf{$\sigma$}_{e}\;  \;\;\forall x \in \Omega^e.
\]
Thus, the spaces $\calU_a $ and $\calS_a$ can be
 numerically discretized to the finite-dimensional spaces $ \calU^{h}_{a}$ and $ \calS^{h}_{a}$, respectively.
In this paper, the shape functions are  based on  piecewise-cubic polynomial for $w(x)$   and piecewise-linear   for $\sig(x)$, i.e.
\[
N_{w}=\left[ \begin{array}{c}\frac{1}{4}\;(1-\xi)^{2}\;(2+\xi)\\ \frac{L_{e}}{8\;}(1-\xi)^{2}\;(1+\xi)  \\  \frac{1}{4}\;(1+\xi)^{2}\;(2-\xi) \\\frac{L_{e}}{8}(1+\xi)^{2}\;(\xi-1) \end{array} \right]   \; ,\;\;\;\;\;\;  N_{\sigma}=\frac{1}{2}\left[ \begin{array}{c}(1-\xi)\\ (1+\xi)\end{array} \right],
\]
where $\xi=2x/L_{e}-1$ with $L_{e}$ is the length of the e-th beam element.
Thus,   on the discretized feasible deformation space $\calU^h_a$, the Gao-Strang total complementary energy can be expressed in the following discretized form
\begin{eqnarray}
 \Xi^h(\bw ,\bsig  ) &=& \sum_{e=1}^{m} \bigg(\frac{1}{2} w^{T}_{e}\; G^{e}(\sigma_e)\; w_{e}-\frac{1}{2} \sigma^{T}_{e}\; K_{e}\; \sigma_e - \lambda^{T}_{e}\; \sigma_{e} -f^{T}_{e}\; w_{e} -c_{e} \bigg)  \nonumber \\
&= &  \frac{1}{2} \bw^{T}\; \bG(\bsig )\; \bw -\frac{1}{2} \bsig^{T}\; \bK\; \bsig  -\blam^{T}\; \bsig  - \bff^{T} \;\bw -c,
 \label{eq:Xi}
 \end{eqnarray}
 where $\bw \in \calU^h_a \subset \calR^{2(m+1)} $ and $\bsig  \in \calS^h_a \subset \calR^{m+1}$ are nodal deflection and dual stress vectors, respectively. We let
\eb
\calS^h_a=\{ \bsig\in \calR^{m+1}| \;\; \det \bG(\bsig) \neq 0  \} .
\ee
The Hessian matrix of the gap function $\bG(\bsig ) \in \calR^{2(m+1)} \times \calR^{2(m+1)}$ is obtained by assembling  the following symmetric matrices $ G^{e}(\sigma_e)$:
\begin{eqnarray}
 G^{e}(\sigma_e) &=& \int_{\Omega_{e}}\bigg( EI\; N''_{w}\; (N''_{w})^{T} + (N_{\sigma})^{T}\; \sigma_{e} \;N'_{w}\; (N'_{w})^{T} \bigg) dx \nonumber \\
&=&  \int^1_{-1}\frac{L_e}{2}\bigg( EI\; N''_{w}\; (N''_{w})^{T} + (N_{\sigma})^{T}\; \sigma_{e} \;N'_{w}\; (N'_{w})^{T} \bigg) d\xi =\bigg[G^e_{ij} \bigg]_{ 4\times 4 },
\end{eqnarray}
where  $G^e_{ij}=G^e_{j i}$  are defined by the two stress ends $\sig^{e}$ and $\sig^{e+1}$ of beam element as:

 \[
 G^e = \left[  \begin{array}{cccc}
  \frac{3(\sig^{e} + \sig^{e+1})}{5Le} + \frac{12EI}{L_e^3}  &
  \frac{\sig^{e+1}}{10}  + \frac{6EI}{L_e^2} &
 -G^e_{11} &
 \frac{\sig^{e}  }{10} + \frac{6EI}{L_e^2} \\
\;\;G^e_{12} & L_e(\frac{\sig^{e}  }{10} + \frac{\sig^{e+1} }{30}) + \frac{4EI}{L_e}  &  - G^e_{1 2} &
  -\frac{L_e}{60} (\sig^{e}  + \sig^{e+1}) + \frac{2EI}{L_e} \\
-G^e_{11} & G^e_{23} &  \;\;G^e_{1 1}  & -G^e_{14} \\
 \;\;G^e_{14}& G^e_{24} &\;\; G^e_{34} &
 L_e(\frac{\sig^{e}  }{30} + \frac{\sig^{e+1} }{10}) + \frac{4EI}{L_e}
\end{array} \right].
\]

The  matrix  $\bK\in\calR^{m+1}\times \calR^{m+1}$ is obtained by assembling the following positive-definite matrices $K_{e} $
\[
 K_{e}=  \int_{\Omega_{e}} \bigg(\frac{3}{2E\alpha}N_{\sigma} \; N^{T}_{\sig} \bigg) dx =  \int^1_{-1} \bigg(\frac{3 L_e}{4E\alpha}N_{\sigma} \; N^{T}_{\sig} \bigg) d\xi  = \frac{L_{e}}{E\alpha} \left[\begin{array}{cc}
 \frac{1}{2} & \frac{1}{4}  \\
\frac{1}{4} &  \frac{1}{2}
 \end{array}\right] .
 \]
Also,  $\blam =  \{\lam_e\} \in \calR^{m+1}$ and $\bff =  \{f_e\}\in \calR^{2(m+1)} $  are defined by assembling the following
  $$  \lambda_e =\int_{\Omega_{e}}\bigg( \frac{3}{2\alpha}\lambda N_{\sig}\bigg) dx =\int^1_{-1}\bigg( \frac{3L_e}{4\alpha}\lambda N_{\sig}\bigg) d\xi = \frac{\lambda L_e}{\alpha}  \left[\begin{array}{cc}
   \frac{3}{4}   \\
\frac{3}{4}
 \end{array}\right],$$
\[
 f_e =  \int_{\Omega_{e}} f\big(x\big)\;  N_w \; dx =\int^1_{-1}\frac{L_e}{2} f\big(\xi\big)\;  N_w d\xi \;,
 \]
 and $c =\sum\limits_{e=1}^m c_e \in \calR$ is defined as
          $$  c_e =\int_{\Omega_{e}} \bigg( \frac{3 E}{4 \alpha} \lambda^2 \bigg) dx=\int^1_{-1}\bigg( \frac{3 E L_e}{8 \alpha} \lambda^2 \bigg) d\xi =  \frac{3 }{4 \alpha} E L_e\lambda^2 . $$

By the critical condition $\delta  \Xi^h(\bw,\bsig)=0$,  canonical equations (\ref{eq-caneqns}) and (\ref{eq-cancons}) have the following discretized forms
\eb\label{eq:ww}
   \bG(\bsig) \; \bw - \bff=0,
\ee
  \eb\label{eq:ww2}
   \frac{1}{2} \bw^T \;\bH  \; \bw - \bK \;\bsig-\blam = 0,
\ee
where $\bH=\bG_{,\bsig}(\bsig)$  stands for gradient of
 $\bG(\bsig)$ with respect to the vector $\bsig$.

For any given  $\bw \in \calU^h_a$, we know that  $\Xi(\bw, *): \calS^h_a \rightarrow \real$ is concave and
 the discretized total potential energy can be obtained by
\eb
\Pi^h_p(\bw) = \max \{ \Xi(\bw, \bsig) | \;\; \bsig\in \calS^h_a \} = \{  \Xi(\bw, \bsig) | \;\;
  \bsig = \bK^{-1} (\half  \bw^T \;\bH  \; \bw  - \blam) \} .
\ee
However, the convexity $\Xi(* , \bsig):\calU_a^h \rightarrow \real$ will depend on $\bsig \in \calS^h_a$.
The discretized pure complementary energy $\Pi^h_d:\calS^h_a \rightarrow \real$ can be obtained by
  the following canonical dual transformation
 \begin{eqnarray}
   \Pi^h_d(\bsig) &=& \mbox{sta } \{ \Xi(\bw, \bsig) | \;\; \bw \in \calU^h_a \} =
   \{ \Xi(\bw, \bsig) | \;\;  \bw= \bG^{-1}(\bsig)  \bff   \} \nonumber \\
 & =&    -\frac{1}{2} \bff^T \; \bG^{-1}(\bsig)\; \bff - \frac{1}{2} \bsig^T \; \bK \; \bsig - \blam^T \; \bsig -c
\label{eq:pi_d}
\end{eqnarray}
where  sta $\{ g(\bw) | \bw \in \calU^h_a\}$ stands for finding the stationary value of $g(\bw)$ on $\calU^h_a$.
 Clearly, its convexity depends on $\bG(\bsig)$. Let
 \eb
   \calS^{+}_a= \{\bsig \in \calS^h_a  \; | \;\bG(\bsig) \succ 0 \},
\ee
 \eb
   \calS^{-}_a = \{\bsig \in \calS^h_a \; | \;\bG(\bsig) \prec 0 \}.  
\ee
Where the symbols ``$\succ$'' and ``$\prec$'' represent to the positive definite matrix and negative definite matrix, respectively.
 \begin{theorem}
Suppose $(\barbw ,\barbsig )$ is a stationary point of $\Xi^h(\bw,\bsig)$, then
$\Pi^h_p(\barbw) = \Xi^h (\barbw, \barbsig) = \Pi^h_d (\barbsig)$.
Moreover, if $\barbsig\in\calS^+_a$, then we have

  {\bf  Canonical Min-Max Duality:}
  $\barbw $ is a global minimizer of $\Pi^h_p(\bw)$ on  $\calU^h_a$  if and only if
  $\barbsig $ is a global maximizer of $\Pi^h_d(\bsig)$ on $\calS^{+}_a$, i.e.,
 \eb\label{eq:Triality1}
 \Pi^h_p(\barbw ) = \min_{\bw\in \calU^h_a} \Pi^h_p(\bw) \;\;\;\Leftrightarrow\;\;\;  \max_{\bsig\in \calS^{+}_a} \Pi^h_d(\bsig)= \Pi^h_d(\barbsig ) . \\
\ee

If $\barbsig \in \calS^{-}_a$, then on a neighborhood $\calU_o \times \calS_o \subset \calU_a^h \times \calS^-_a$ of $(\barbw, \barbsig)$ we have

  {\bf Canonical Double-max Duality:}
  The stationary point $\barbw $ is a local maximizer of $\Pi^h_p(\bw)$ on  $\calU_o$ if and only if the stationary point $\barbsig $ is a local maximizer of $\Pi^h_d(\bsig)$ on $ \calS_o$, i.e.,
 \eb\label{eq:Triality2}
 \Pi^h_p(\barbw ) = \max_{\bw\in \calU_o} \Pi^h_p(\bw) \;\;\;\Leftrightarrow\;\;\;  \max_{\bsig\in \calS_o} \Pi^h_d(\bsig)= \Pi^h_d(\barbsig )  \\
\ee

 {\bf Canonical Double-min Duality (if $\mbox{dim} \calU^h_a  = \mbox{dim}  \calS^h_a $):}
  The stationary point $\barbw $ is a local minimizer of $\Pi^h_p(\bw)$ on  $\calU_o$ if and only if the stationary point $\barbsig $ is a local minimizer of $\Pi^h_d(\bsig)$ on $ \calS_o$, i.e.,
 \eb\label{eq:Triality2}
 \Pi^h_p(\barbw ) = \min_{\bw\in \calU_o} \Pi^h_p(\bw) \;\;\;\Leftrightarrow\;\;\;  \min_{\bsig\in \calS_o} \Pi^h_d(\bsig)= \Pi^h_d(\barbsig )
\ee
\end{theorem}

The proof of this theorem follows from the general results in global optimization
\cite{chen-gao-jogo,gao-wu-jimo,mora-gao-memo}.
Canonical min-max duality
 can be used to find  the global minimizer of the nonconvex problem via the following canonical dual problem:
  \eb
   (\calP^d): \;\;\; \max \{ \Pi^h_d(\bsig) | \; \bsig \in \calS^{+}_a  \} , \label{eq-cdmax}
   \ee
which is a concave maximization problem and can be solved easily by well-developed convex analysis and optimization techniques.
    The   canonical double-max  and  double-min duality statements   can be used to find
     the biggest local maximizer and a local minimizer of the nonconvex  primal problem, respectively.
It was proved in \cite{chen-gao-jogo,gao-wu-jimo,mora-gao-memo} that
both  canonical min-max and double-max duality statements hold strongly regardless of the dimensions of
$\calU^h_a$ and $\calS^h_a$, while the canonical double-min duality statement (\ref{eq:Triality2}) holds strongly  for
$\dim \calU^h_a  = \dim \calS^h_a$, but weakly
if  $\dim \calU^h_a  \neq \dim \calS^h_a$. This case is within our reach in the following applications.

\section{ Semi-Definite Programming Algorithm}
\label{sec:Applications}  
It is easy to understand that the nonconvex post-buckling problem could have multiple global minimizers for certain
external loads, say $q(x)= 0$. In this case  we have $\det \bG(\bsig) = 0$ and $\calS_a^h = \emptyset$.
In order to deal with this case, this section presents a SDP (Semi-Definite Programming, see  \cite{Bernd-Matousek} and \cite{Stephen-Lieven}) reformulation to solve the canonical dual problem  (\ref{eq-cdmax}). The SDP algorithm is applied to obtain  all post-buckled solutions of a large deformed elastic  beam. 

By the fact that $ \Xi(\bw,\bsig)$ is a saddle function on $\calU_a^h \times \calS^+_a$,
 we have
\eb
\min_{\bw \in \calU^h_a} \Pi^h_p (\bw) = \min_{\bw \in \calU^h_a} \max_{\bsig \in \calS^+_a} \Xi(\bw, \bsig) =
 \max_{\bsig \in \calS^+_a} \min_{\bw \in \calU^h_a}  \Xi(\bw, \bsig) .
\ee
For any given $\bsig \in \calS^+_a$, the solution to $ \min_{\bw \in \calU^h_a}  \Xi(\bw, \bsig) $ leads to
\eb  \label{eq:w(sig)}
\bw=\bw(\bsig)= \bG^{-1}(\bsig)\bff
\ee
Thus, the stress fields $\bsig$
can be found by the following problem
\begin{eqnarray}
  \max_{\bsig }\;\Xi(\bw(\bsig), \bsig)
 &=&\frac{1}{2} \bw(\bsig)^T \bG(\bsig) \bw(\bsig)   - \frac{1}{2} \bsig^T \; \bK \; \bsig -\blam^T \bsig -\bff^T \bw(\bsig) -c    \nonumber \\
 &\equiv &  \max_{\bsig }\; \Pi^h_d(\bsig)   \nonumber \\
\mbox{ s.t. } \;\;\; \;\;\;\;\;\;\;\; \bG(\bsig) &\succeq& 0,
  \label{eq:Xi(w(sig))}
 \end{eqnarray}
 where the symbol ``$\succeq$'' represents to the positive semi-definite matrix.
  By   canonical min-max duality we know that if
   $\bsig^* \in \calS^{+}_a$ is a global maximizer of  problem (\ref{eq:Xi(w(sig))}),
    then  $\bw^*=\bw(\bsig^*)$ should be a global minimizer of $\Pi^h_p(\bw)$.
Furthermore, the problem (\ref{eq:Xi(w(sig))})  is the same as:
\begin{eqnarray}
  & &\max_{\bsig, t }\;t   \;\; \mbox{ s.t. } \;\;  \bG(\bsig) \succeq 0,
 \;\;\;\;\;  t \leq  \phi( \bsig) - \frac{1}{2} \bsig^T \; \bK \; \bsig
  \label{eq:t}
 \end{eqnarray}
 where $\phi( \bsig) = \frac{1}{2} \bw(\bsig)^T \bG(\bsig) \bw(\bsig) -  \blam^T \bsig -\bff^T \bw(\bsig) -c$.
 By the fact that $\bK\succ 0$,  the Schur complement lemma (see \cite{Stephen-Lieven})  for  the second inequality constraint in
 (\ref{eq:t}) implies
 \begin{eqnarray}
  \left[\begin{array}{cc}
2 \bK^{-1} &  \bsig  \\
\bsig^T &
\phi( \bsig) - t
 \end{array}\right] \succeq 0  .\label{eq:GlobalMax}
\end{eqnarray}
Thus,  the problem (\ref{eq:t})  can be relaxed to the following   Semi-Definite Programming (SDP)  problem
\begin{eqnarray}
& & \max_{ \bsig, t}   \; t   \;\; \; \mbox{  s.t. } \;\;  \bG(\bsig) \succeq 0, \;\; \;
  \left[\begin{array}{cc}
 2 \bK^{-1} &  \bsig  \\
\bsig^T & \phi(  \bsig) - t
 \end{array}\right] \succeq 0 . \label{eq:GlobalMax}
\end{eqnarray}

In the same way, the SDP relaxation for the canonical double-max duality statement
\eb
\max_{\bw \in \calU^h_a} \Pi^h_p(\bw) = \max_{\bw , \bsig} \Xi(\bw, \bsig) = \max  \Pi^h_d(\bsig)  \;\;
\mbox{ s.t. } \;\;\bsig \in \calS^-_a
\ee
  should be equivalent to
  \begin{eqnarray}
& & \max_{\bsig , t}   \; t   \;\; \mbox{  s.t. } \;\;  - \bG(\bsig) \succ  0, \;\; \;
  \left[\begin{array}{cc}
2  \bK^{-1} &  \bsig  \\
\bsig^T &
 \phi( \bsig)  - t
 \end{array}\right] \succeq 0  , \label{eq:LocalMax}
\end{eqnarray}
 which leads to a local maximum solution to the post-buckling problem.

Now, let $(\bw^*,\bsig^*)$ be  a  local minimizer of the canonical double-min  problem
 $\; \min_{\bw } \Pi^h_p(\bw) = \min_{\bw} \max_{\bsig} \Xi(\bw, \bsig) = \min_{\bsig} \Pi^h_d(\bsig) $ s.t. $\bsig \in \calS^-_a$.
By eq.(\ref{eq:w(sig)}), the local minimizer is  equivalent to the following problem
\begin{eqnarray}
  & &\min_{\bsig }\;\{\Xi(\bw(\bsig), \bsig) \equiv \Pi^h_d(\bsig)\}   \;\;\; \mbox{ s.t. } \;\;
   \;  \bG(\bsig) \prec 0.
 \end{eqnarray}
 This problem is  the same as:
\begin{eqnarray}
& & \min_{\bsig, t}   \; t   \;\; \mbox{  s.t.  } \;\;  \bG(\bsig)  \prec  0, \;\; \;
t \ge -\frac{1}{2} \bff^T \; \bG^{-1}(\bsig)\; \bff - \frac{1}{2} \bsig^T \; \bK \; \bsig - \blam^T \; \bsig -c .  \label{minXi22}
\end{eqnarray}
  In order to  apply the  Schur complement lemma for the second inequality in  (\ref{minXi22}),
we  need to  linearize the complementary energy   $V^*(\bsig) = \frac{1}{2} \bsig^T \; \bK \; \bsig $.
This can be done by  using a  reformulated  pure complementary energy:
 \eb\label{eq:New_pi_d}
 \widehat{\Pi}^d(\bsig,\bw)= -\frac{1}{2} \bff^T \; \bG^{-1}(\bsig)\; \bff - \frac{1}{2} \bw^T \; \bM(\bsig) \; \bw -\frac{1}{2} \blam^T \; \bsig -c ,
\ee
the stiffness matrix  $\bM(\bsig)$ in the strain energy $V(\bw) = \frac{1}{2} \bw^T \; \bM(\bsig) \; \bw = V^*(\bsig)$
is obtained by assembling the following symmetric matrices $ M^e(\sig_e)$ in each element
\begin{eqnarray}
      M^e(\sig_e) &=& \int_{\Omega_{e}} \frac{1}{2} \bigg(  (N_{\sigma})^{T}\; \sigma_{e} \;N'_{w}\; (N'_{w})^{T} \bigg) dx  = \int^1_{-1}\frac{L_e}{4}\bigg( (N_{\sigma})_{T}\; \sigma_{e} \;N'_{w}\; (N'_{w})^{T} \bigg) d\xi \nonumber \\
  &=&  \left[  \begin{array}{cccc}
 \frac{3}{10Le} (\sig^{e}  + \sig^{e+1}) &
 \frac{1}{20}\sig^{e+1}& \; - M^e_{1,1} &
\frac{1}{20}\sig^{e} \\
 M^e_{12} &  \frac{L_e}{60} (3\sig^{e} + \sig^{e+1})  &  - M^e_{1 2} &
   \frac{-L_e}{120} (\sig^{e}  + \sig^{e+1}) \\
M^e_{13} & M^e_{23} &  \;\;M^e_{1 1}  & -M^e_{14} \\
 M^e_{14}& M^e_{24} & \;\;M^e_{34} &
  \frac{L_e}{60}(\sig^{e}  + 3\sig^{e+1})
\end{array} \right].\label{eq:M}
\end{eqnarray}
 Therefore,  by using   $\widehat{\Pi}^d(\bsig,\bw)$,
  problem (\ref{minXi22})  can be relaxed to
\begin{eqnarray}
& & \min_{\bsig, t}   \; t   \;\; \mbox{  s.t. } \;\;   \bG(\bsig)  \prec  0, \;\; \;
\frac{1}{2} \bff^T \; \bG^{-1}(\bsig)\; \bff + \hat{\phi} (\bsig,\bw) +t \ge 0,
\label{newPure}
\end{eqnarray}
where $\hat{\phi}(\bsig,\bw) = \frac{1}{2} \bw^T \; \bM(\bsig) \; \bw +\frac{1}{2} \blam^T \; \bsig +c $.
The primal variable $\bw$ in this problem can be computed by the dual solution $\bsig$ in the primal-dual iteration.
Thus, by using  the Schur complement lemma  this problem can be relaxed to the following   SDP  problem
\begin{eqnarray}
& & \min_{ \bsig, t}   \; t   \;\;  \mbox{ s.t.  } \;\;   - \bG(\bsig) \succ  0, \;\; \;
  \left[\begin{array}{cc}
-2 \bG(\bsig) &  \bff  \\
\bff ^T &    \hat{\phi}(\bsig, \bw)  +t
 \end{array}\right] \succeq 0 , \label{eq:LocalMin}
\end{eqnarray}
 Clearly, if  stress  $\bsig^* $ is a local minimizer on $\calS^{-}_a$ of   problem (\ref{eq:LocalMin}),
  the canonical double-min duality shows that   $\bw^*=\bw(\bsig^*)$ should be  a local minimizer of $\Pi^h_p(\bw)$.

Consequently, the primal-dual semi-definite programming (PD-SDP) algorithm
  for solving all possible  post-buckling solutions can be proposed as the following.  \vspace{.2cm} \\
  {\bf PD-SDP Algorithm}:
\begin{enumerate}
  \item Given  initial primal solution  $\bw^{(0)}$ and error allowance   $ \epsilon > 0$ . Let $k=1$ ;

    \item Compute the dual solutions $\{\bsig^{(k)} \}$  by applying the SDP solver to problems
      (\ref{eq:GlobalMax}),   (\ref{eq:LocalMax}) and
  (\ref{eq:LocalMin}), respectively.

  \item Compute the primal solution $\bw^{(k)}= [\bG(\bsig^{(k)})]^{-1}\bff$.
  \item For check convergence;
   if $ \| \bw^{(k)}- \bw^{(k-1)}  \| /  \|\bw^{(k)} \| \leq \epsilon$, stop with the optimal solution $\bw^* =\bw^{(k)} $. Otherwise, let  $k=k+1$ and go to step 2.
\end{enumerate}

The SDP solver  used in this algorithm is  a popular software package named SeDuMi, which is
 based on the  interior point method \cite{Sturm-99}.

\section{Numerical solutions} \label{sec:Numerical solutions}
We present in this section two different types of  beams. Geometrical data were kept fixed for all computations; elastic modulus $E=1000 Pa$,   Poisson's  ratio $\mu=0.3$ and beam length $L=1 \mathrm{m}$.  The lateral load $q(x)$ is assumed to be either a uniformly distributed load such that $f(x)=(1-\mu^2) q(x)=0.1 N/\mathrm{m}$ or a concentrated force on the center of the beam in which $f(x)=0.1 N$.
A different numbers of elements with the same beam length,  
different compressive load $\lambda$  with different values of beam height 
are applied in this paper.
 \begin{figure}[h!]
\begin{center}
 \scalebox{0.18}
{\includegraphics {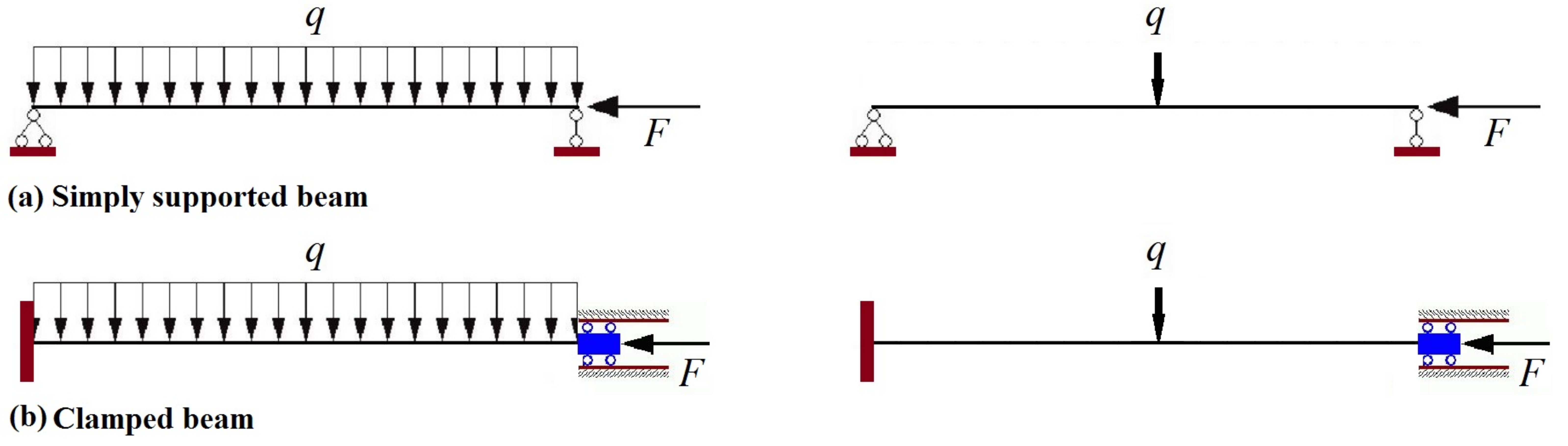}}
\caption{{\em Types of beams - uniformly distributed load (left)}, {\em concentrated force (right)}}
\label{SSB}
\end{center}
\end{figure}
\vspace{-0.5cm}
\subsection{Simply supported beam}
A simply supported beam model is fixed in both directions at $x=0$ and fixed only in the y-direction at $x=L$ as shown in
Figure (\ref{SSB}-a) with the boundary conditions $w(0)=w''(0)=w(L)=w''(L)=0$. 
 \begin{figure}[h!]
\begin{center}
\scalebox{0.155}{\includegraphics{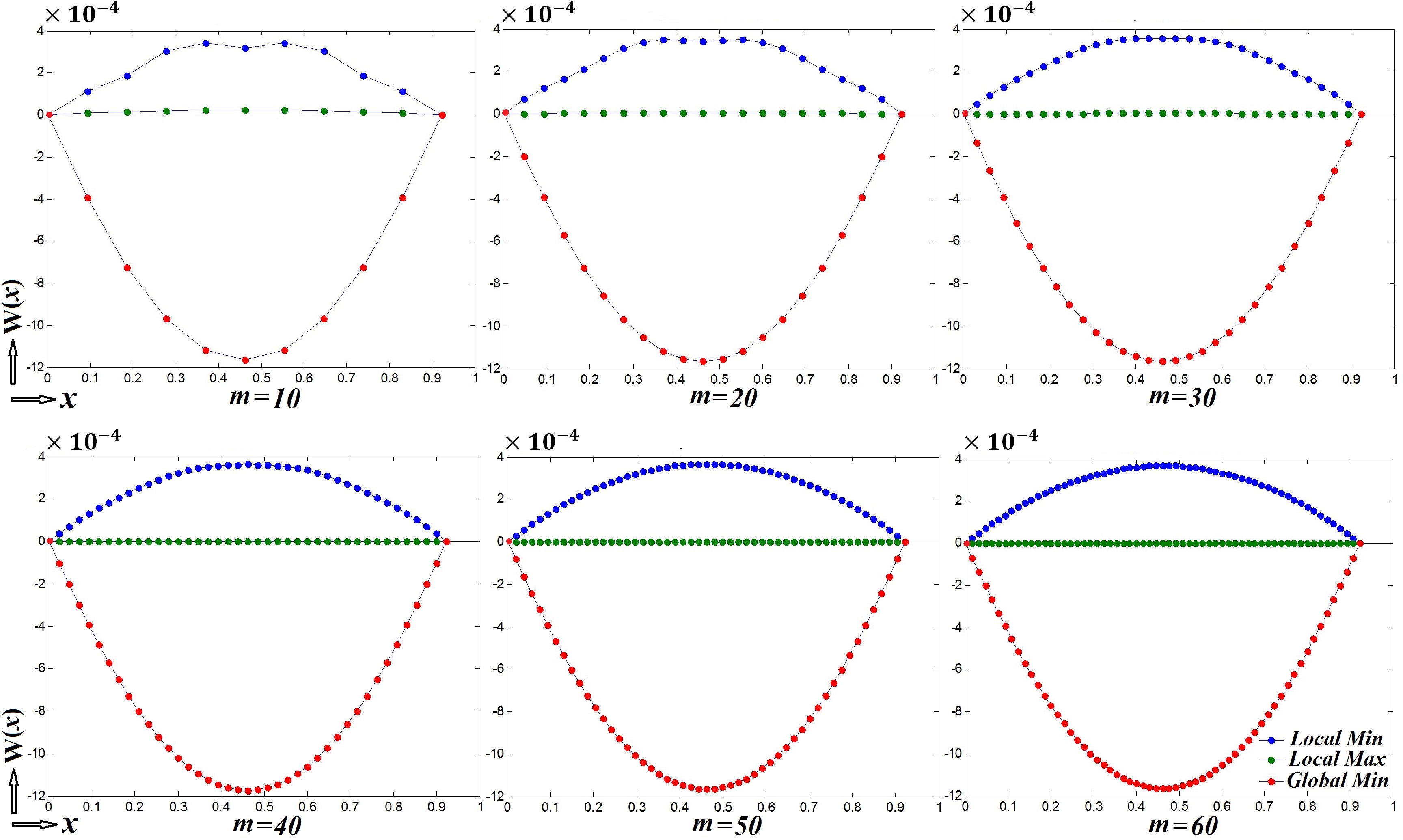}}
\caption{{\em \footnotesize{Simply supported beam under a uniformly distributed load with $\lambda=0.01\mathrm{m}^2$ ($h=0.05\mathrm{m}$)}}}
\label{a1}
\end{center}
\end{figure}
If the beam height is 0.1 (i.e. $h=0.05\mathrm{m}$), the critical load  is  $\lambda_{cr}= 0.00097\mathrm{m}^2$  (see eq.(\ref{eq:Fcr})). 
For a different numbers of beam elements, the approximate deflections of this beam  with $\lambda=0.01\mathrm{m}^2\textgreater \lambda_{cr}$ under a uniformly distributed load   are illustrated in Figure \ref{a1}.
In the graphs, 
red represents the global minimum,
green represents the local maximum and
blue represents the local minimum of $\Pi(w)$. 
Figure \ref{a1} shows that the two post-buckled configurations; global minimum and local maximum, look alike with all of the different numbers of beam elements. 
In contrast to the local minimum,   few  differences appear in the local unstable buckled configuration. 
The curve charts with 40, 50 and 60 elements seem very similar and more stable than the curve charts that contain 10, 20 and 30 elements.
 Once again, Figure \ref{a2}  shows that, 
with a different number of elements at $\lambda=0.015\mathrm{m}^2\textgreater \lambda_{cr}$,
slight differences appear on the local minimum curves.
  \begin{figure}[h!] 
\begin{center}
\scalebox{0.155}{\includegraphics{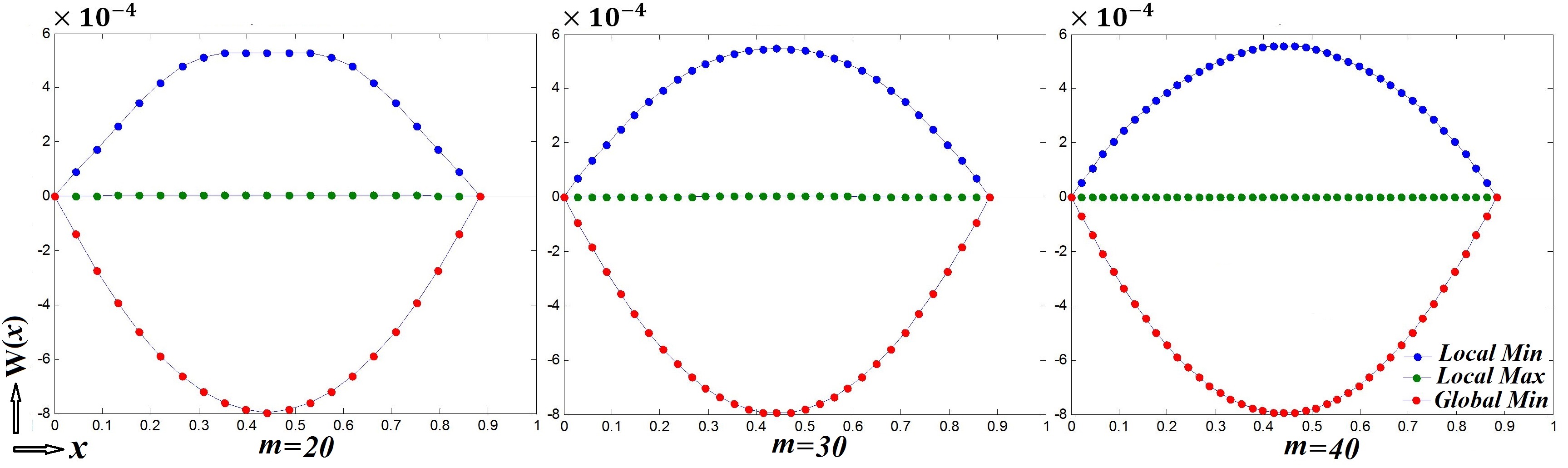}}
\caption{{\em \footnotesize{Simply supported beam under a uniformly distributed load with} $\lambda=0.015\mathrm{m}^2$ ($h=0.05\mathrm{m}$)}}
\label{a2}
  \end{center}
  \end{figure}
  
\begin{figure}[h!]
  \begin{center}
\scalebox{0.156}{\includegraphics{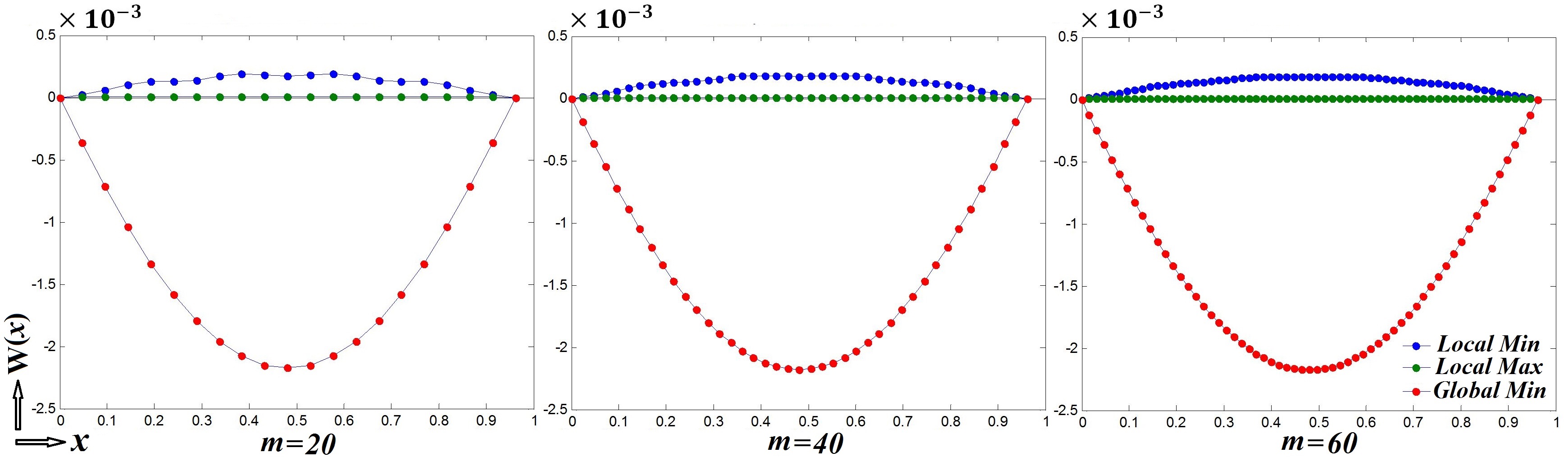}}
\caption{ {\em \footnotesize{Simply supported beam under a uniformly distributed load with $\lambda=0.005\mathrm{m}^2$ ($h=0.05\mathrm{m}$)}}}
\label{a3}
  \end{center}
  \end{figure}

   The local minimum solutions with a different number of beam elements
at a compressive load $\lambda=0.005\mathrm{m}^2\textgreater \lambda_{cr}$
look alike,  as shown in Figure  \ref{a3}.
     The Gao-Strang gap function
    for all post-buckled solutions was computed under a uniformly distributed load for a different number of elements  with $\lambda=0.01\mathrm{m}^2$ as reported in Table \ref{SSB00}. 

\begin{figure}[h!]
  \begin{center}
\scalebox{0.115}{\includegraphics{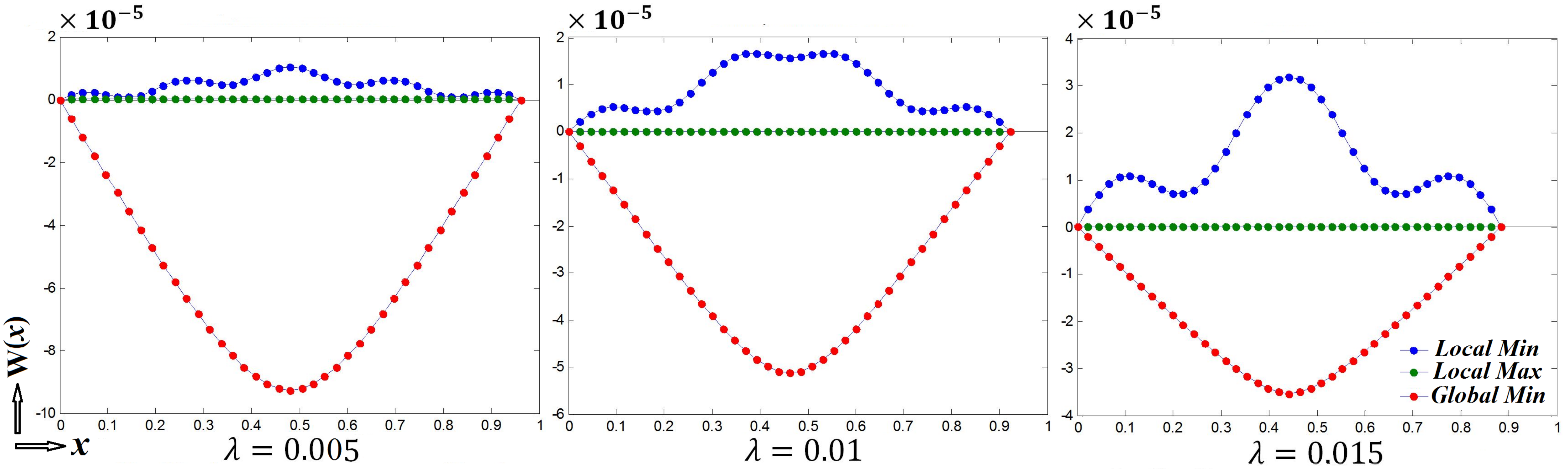}}
\caption{{\em \footnotesize{Simply supported beam under a concentrated force  ($h=0.05\mathrm{m}$)}}}
\label{a4}
  \end{center}
  \end{figure}

\begin{figure}[h!]
  \begin{center}
\scalebox{0.115}{\includegraphics{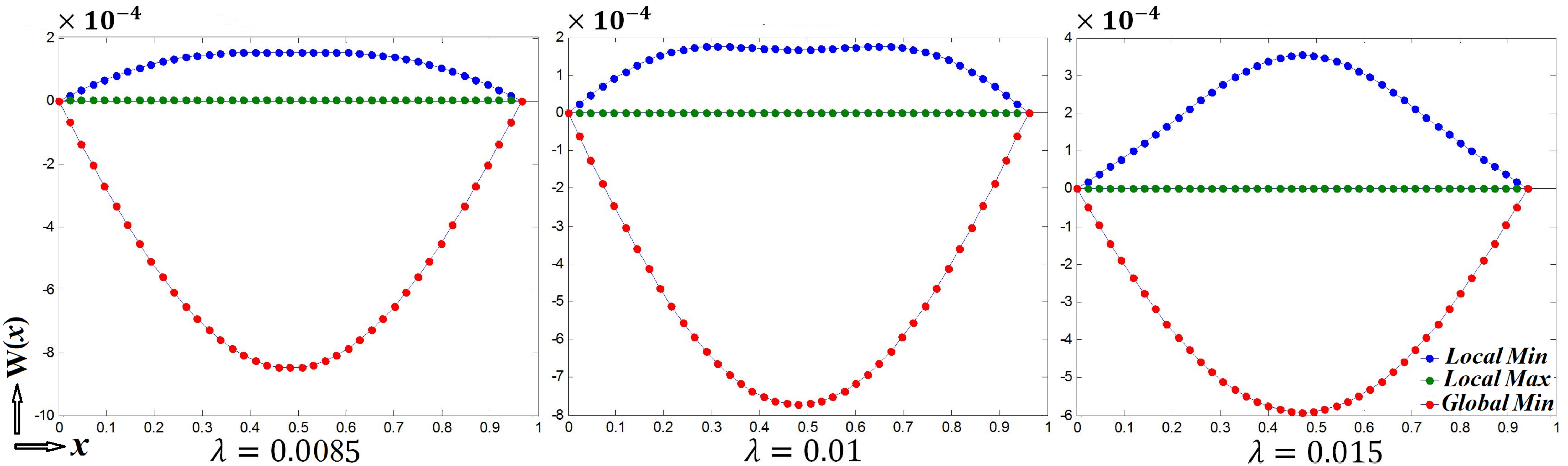}}
\caption{{\em \footnotesize{Simply supported beam under a uniformly distributed load  ($h=0.1\mathrm{m}$)}}}
\label{a5}
  \end{center}
  \end{figure}

\begin{figure}[h!]
  \begin{center}
\scalebox{0.115}{\includegraphics{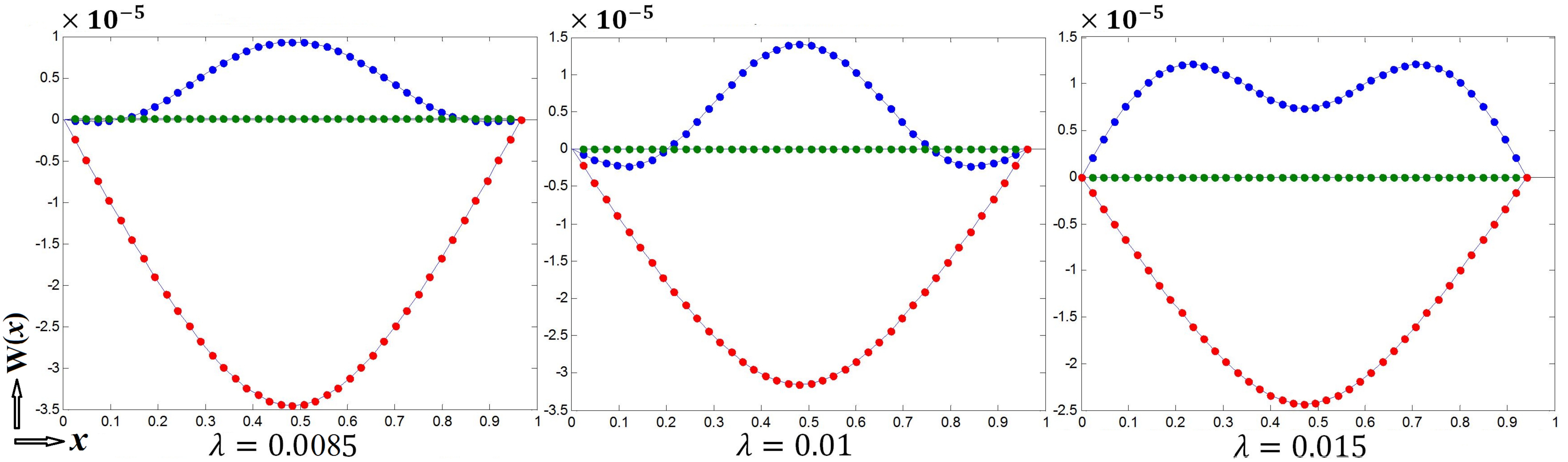}}
\caption{{\em \footnotesize{Simply supported beam with  a concentrated force ($h=0.1\mathrm{m}$)}}}
\label{a6}
\end{center}
\end{figure}

We  focus on 40 elements with the same beam length for all the following examples.
The   deflections of the simply supported beam under  a concentrated force with different
 compressive loads $\lambda  >  \lambda_{cr}$ are  illustrated in  Figure \ref{a4}.
  At  $h=0.1\mathrm{m}$, the critical load of the simply supported beam is  $\lambda_{cr}= 0.0078\mathrm{m}^2$. The  deflections of this beam under a uniformly distributed load and a concentrated force  are summarized in Figures \ref{a5}  and  \ref{a6},  respectively.
       The Gao-Strang gap function
    for all three post-buckled solutions  was computed under a uniformly distributed load and a concentrated force 
 as reported in Tables \ref{SSB1} and \ref{SSB2}, respectively.
\begin{table}[h!]
\begin{center}
\begin{tabular}{ |p{2cm}||p{1.5cm}|p{2.5cm}|p{2.5cm}|p{2.5cm}|}
\hline
    \multicolumn{1}{|c||}{\small{Compressive}}
    & \multicolumn{1}{|c|}{\small{No elements}}& \multicolumn{3}{|c|}{\small{Gap function under a uniformly distributed load}}\\
 \cline{3-5}
  \small{load}  & \centerline{$m$} & \small{Global Min} &   \small{Local Min} & \small{Local Max}\\
 \hline 
                           &   $\centerline{20}$    &  7.63568e-09     & -2.15332e-09    &   -4.16926e-07    \\
 $\lambda=0.01$  &   $\centerline{40}$     &  1.45323e-09     & -8.56515e-10    &  -1.04182e-07     \\
                          &  $\centerline{60}$   &  6.10785e-10   &  -4.93895e-10   &       -4.62995e-08  \\
  \hline 
\end{tabular}
        \caption{{\em \footnotesize{Gao-Strang gap function for simply supported beam with different numbers of elements}\label{SSB00}}}
        \end{center}
\end{table}
\begin{table}[h!]
\begin{center}
\begin{tabular}{ |p{2cm}||p{1.5cm}|p{2.5cm}|p{2.5cm}|p{2.5cm}|}
\hline
    \multicolumn{1}{|c||}{}
    & \multicolumn{1}{|c|}{\small{Compressive }}& \multicolumn{3}{|c|}{\small{Gap function under a uniformly distributed load}}\\
 \cline{3-5}
  \small{Beam height}  &\small{ loads ``$\lambda$''} & \small{Global Min} &   \small{Local Min} & \small{Local Max}\\
 \hline 
 
                 &   0.005     &  1.38767e-09     & -3.90449e-10    &   -1.04182e-07    \\
 $h=0.05$   &   0.01     &  1.45323e-09      & -8.56515e-10     &  -1.04182e-07      \\
                 &   0.015    &  1.51964e-09   &  -1.01164e-09   &       -1.04182e-07   \\
  \hline 
                      &   0.0085    &  1.66885e-10    &  -1.48050e-10    &  -1.30228e-08    \\
  $h=0.1$      &   0.01     &  1.67195e-10    &  -1.50613e-10     &       -1.30228e-08    \\
                       &   0.015     &  1.68227e-10   & -1.55455e-10    &   -1.30228e-08   \\
 \hline
\end{tabular}
        \caption{{\em \footnotesize{Gao-Strang gap function for simply supported beam under a uniformly distributed load}\label{SSB1}}}
        \end{center}
\end{table}
\begin{table}[h!]
\begin{center}
\begin{tabular}{ |p{2cm}||p{1.5cm}|p{2.5cm}|p{2.5cm}|p{2.5cm}|}
\hline
    \multicolumn{1}{|c||}{}
    & \multicolumn{1}{|c|}{\small{Compressive}}& \multicolumn{3}{|c|}{\small{Gap function under a concentrated load}}\\
 \cline{3-5}
  \small{Beam height}  & \small{ loads ``$\lambda$''}& \small{Global Min} &   \small{Local Min} & \small{Local Max}\\
 \hline 
 
                 &   0.005     &   2.72407e-12  & -9.56982e-13   &-1.89005e-10      \\
 $h=0.05$   &   0.01     &  2.84230e-12   &  -1.78093e-12   &   -1.89005e-10   \\
                 &   0.015    &  2.96381e-12   & -2.05556e-12   &  -1.89005e-10   \\
  \hline 
                      &   0.0085    & 3.28941e-13    & -2.95470e-13   &   -2.36257-11  \\
  $h=0.1$      &   0.01     &  3.29501e-13   & -3.00014e-13   &   -2.36257e-11   \\
                       &   0.015     &  3.31372e-13   & -3.08597e-13   & -2.36257e-11   \\
 \hline
\end{tabular}
        \caption{{\em \footnotesize{Gao-Strang gap function for simply supported beam under a concentrated load}\label{SSB2}}}
        \end{center}
\end{table}

 \subsection{Doubly/Clamped beam}
\begin{figure}[h!]
  \begin{center}
\scalebox{0.16}{\includegraphics{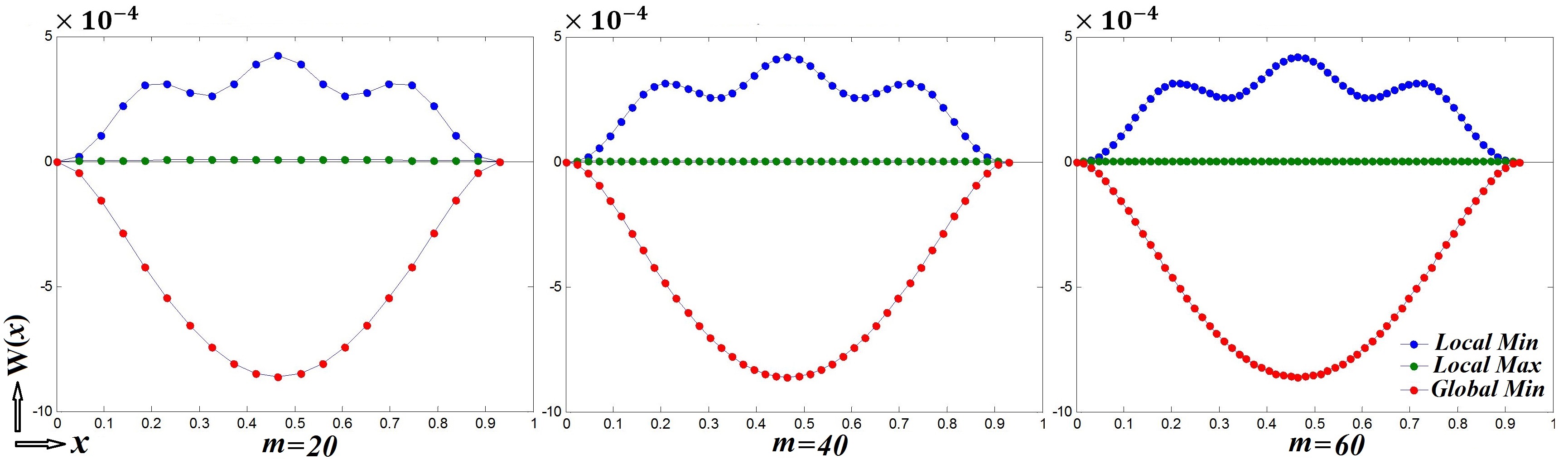}}
\caption{{\em \footnotesize{Clamped beam under a uniformly distributed load with $\lambda=0.009\mathrm{m}^2$ ($h=0.05\mathrm{m}$)}}}
\label{U0}
\scalebox{0.120}{\includegraphics{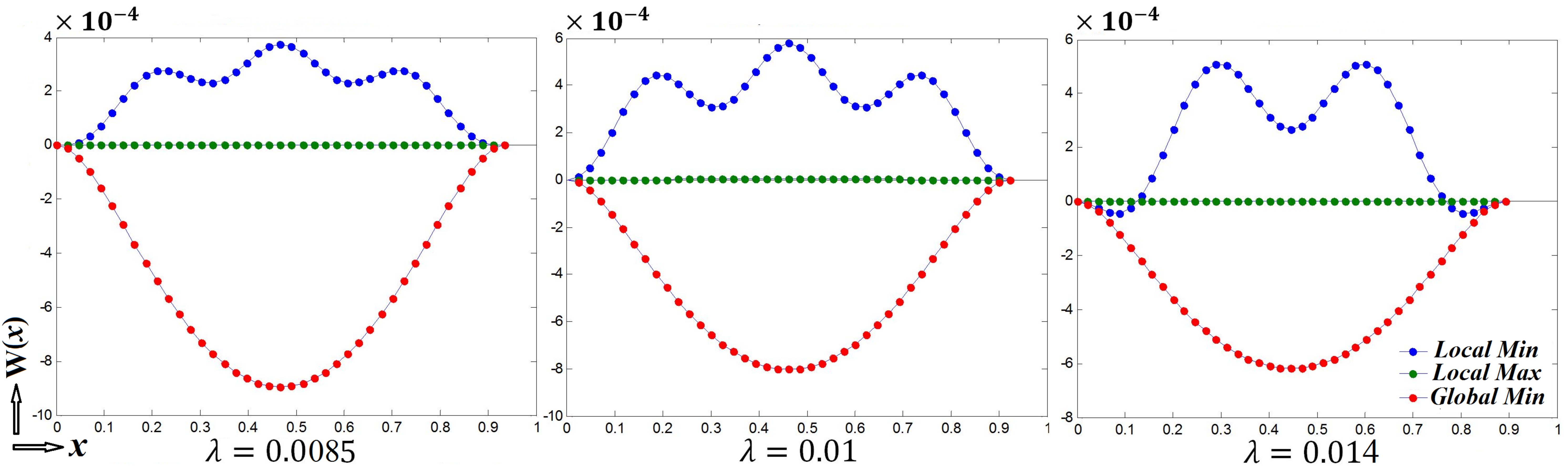}}
\caption{{\em \footnotesize{Clamped beam under a uniformly distributed load ($h=0.05\mathrm{m}$)}}}
\label{b1}
 \scalebox{0.120}{\includegraphics{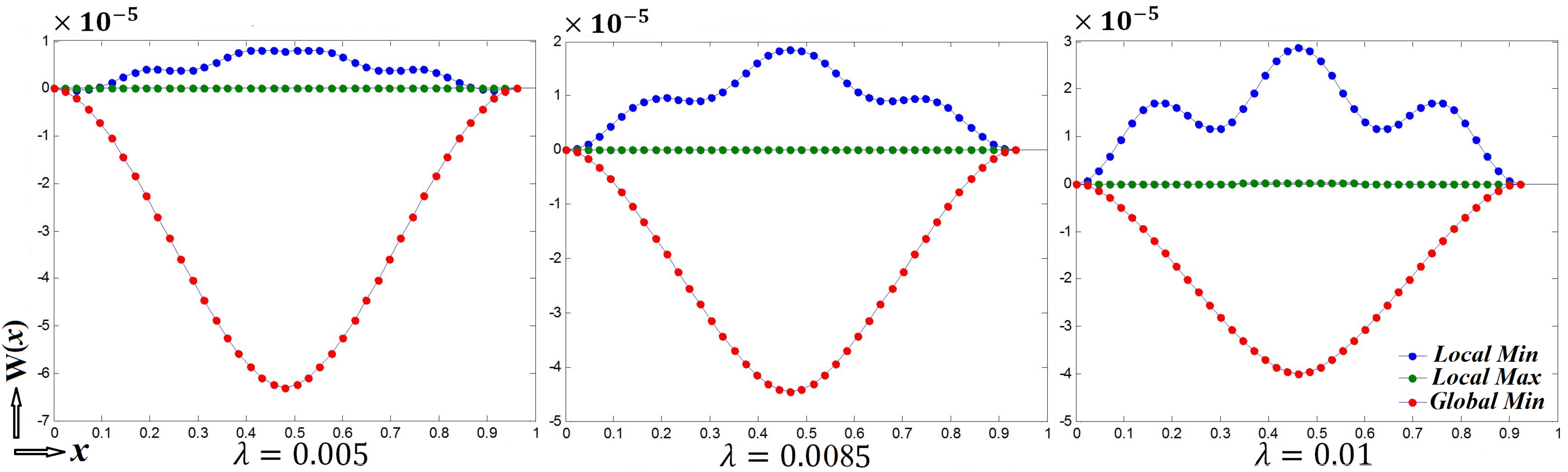}}
 \caption{{\em \footnotesize{Clamped beam under a concentrated force ($h=0.05\mathrm{m}$)}}}
 \label{b2}
\end{center}
\end{figure}

A clamped beam or doubly/clamped beam model  is clamped  at both ends as shown in Figure (\ref{SSB}-b). The boundary conditions are defined as; $w(0)=w'(0)=w(L)=w'(L)=0$.
The  Euler buckling   load  of this beam with  $h=0.05\mathrm{m}$ is  $\lambda_{cr}= 0.0041\mathrm{m}^2$.
A different number of beam elements are applied  with the same conditions and $\lambda=0.009\mathrm{m}^2$. We found  that  the results looked alike for all three post-buckled solutions as shown in Figure \ref{U0}.
   The results of the deflections under a uniformly distributed load and a concentrated force for   different axial loads $\lambda \textgreater \lambda_{cr}$ with $m=40$ are illustrated in Figures \ref{b1}  and  \ref{b2}, respectively.     

        The Gao-Strang gap function
    for all three post-buckled solutions,  with  different axial loads and beam heights, was computed
     under a uniformly distributed load and a concentrated force 
      as reported in Tables \ref{BC41} and \ref{BC42}, respectively. 

\begin{table}[h!]
\begin{center}
\begin{tabular}{ |p{2cm}||p{1.5cm}|p{2.5cm}|p{2.5cm}|p{2.5cm}|}
\hline
    \multicolumn{1}{|c||}{}
    & \multicolumn{1}{|c|}{\small{Compressive}}& \multicolumn{3}{|c|}{\small{Gap function under a uniformly distributed load}}\\
 \cline{3-5}
  \small{Beam height}  & \small{ loads ``$\lambda$''} & \small{Global Min} &   \small{Local Min} & \small{Local Max}\\
 \hline 
 
                 &   0.0085    &  2.09541e-08     & -2.01747e-08    &-1.04101e-07     \\
 $h=0.05$   &  0.009     &   2.09619e-08    &  -2.02106e-08   &  -1.04101e-07   \\
                 &    0.01   &  2.09768e-08  &  -2.02717e-08   &    -1.04101e-07     \\
                 &    0.014   & 2.10396e-08   &   -2.04287e-08  &   -1.04101e-07      \\
 \hline
\end{tabular}
        \caption{{\em \footnotesize{Gao-Strang gap function for  doubly/clamped beam under a uniformly distributed load}\label{BC41}}}
        \end{center}
\end{table}
\begin{table}[h!]
\begin{center}
\begin{tabular}{ |p{2cm}||p{1.5cm}|p{2.5cm}|p{2.5cm}|p{2.5cm}|}
\hline
    \multicolumn{1}{|c||}{}
    & \multicolumn{1}{|c|}{\small{Compressive}}& \multicolumn{3}{|c|}{\small{Gap function under a concentrated load}}\\
 \cline{3-5}
  \small{Beam height}  & \small{ loads ``$\lambda$''}& \small{Global Min} &   \small{Local Min} & \small{Local Max}\\
 \hline 
 
                 &   0.005   & 1.08569e-11      &  -9.01280e-12   &    -1.88954e-10 \\
 $h=0.05$   &  0.0085     &  1.09445e-11     & -9.72096e-12    &   -1.88954e-10  \\
                 &    0.01   &  1.09801e-11  & -9.87268e-12    &  -1.88954e-10       \\
 \hline
\end{tabular}
        \caption{{\em \footnotesize{Gao-Strang gap function for  doubly/clamped beam under a concentrated load}\label{BC42}}}
        \end{center}
\end{table}
\section{Conclusions}
We have presented a canonical dual finite element method for the
post-buckling analysis of a  large  deformed elastic  beam proposed by Gao in 1996.
The nonconvexity of the  total potential energy $\Pi(w)$   is necessary for  the post-buckling phenomenon,
but it leads to a fundamental difficulty for traditional numerical methods and algorithms.
Based on the canonical duality theory and mixed finite element method, a new primal-dual semi-definite program (PD-SDP)
algorithm is proposed, which can be used  to
solve this   challenging nonconvex variational problem to obtain all
possible post-buckled solutions.
 Extensive applications are illustrated for the post-buckled beam with
   different boundary conditions and axial compressive forces.
   The Gao-Strang gap function is computed for all post-buckled solutions. It is interesting to note that for local and global minima,  the value of this gap function is affected by  both the number of beam elements and axial loads, but 
for local maxima, its value is affected mainly  by the number of elements.   
Our results show that  the number of post-buckling solutions depends mainly on the axial compressive forces.
      For a given nontrivial $q(x)$, the nonlinear beam can have at most three post-buckled solutions if  $\lambda \ge \lambda_{cr}$.
 Both the global minimizer and local maximizer solutions are very stable.
  However, the local minimal solution is very sensitive not only to the artificial parameters, such as  the size  of the finite elements,  
   but also to  the natural conditions such as  the axial compressive forces and  boundary conditions. 
    Particularly, for a given $\lambda > \lambda_{cr}$, the biger is the external load $q(x)$, the smaller is the local minimal solution $w(x)$.
Therefore, the related  numerical results presented  in  Figure 13 in \cite{cai-gao-qin} are wrong\footnote{The  first author of \cite{cai-gao-qin} is responsible for this mistake since he didn't let the other two co-authors to check his computer code.}. 
\\

{\bf Acknowledgements}: The authors would like to sincerely acknowledge  the important comments and suggestions  from an anonymous reviewer, which significantly improved the quality of the manuscript. 
This  research was supported by the US Air Force Office of Scientific Research under
 the grants  (AOARD) FA2386-16-1-4082 and  FA9550-17-1-0151.   

\end{document}